\documentclass[10pt]{amsart}
\usepackage{amssymb,amsmath,amsfonts,amsthm,amscd,enumerate}

\usepackage{graphics}

\input xy
\xyoption{all}

\newcommand{\CC}{\mathbb C}
\newcommand{\RR}{\mathbb R}
\newcommand{\maA}{\mathcal{A}}

\newcommand{\maD}{\mathcal{D}}

\newcommand{\ZZ}{\mathbb{Z}}

\newcommand{\HH}{\mathbb{H}}

\newcommand{\smooth}{\mathcal{C}^{\infty}}

\theoremstyle{definition}
\newtheorem{definition}{Definition}[section]
\newtheorem{exm}[definition]{\bf Example}

\theoremstyle{plain}

\newtheorem{prop}[definition]{Proposition}
\newtheorem{lem}[definition]{Lemma}

\theoremstyle{remark}

\begin{document}
\title{Notes on generalized pseudo-differential operators}

\author{Shantanu Dave}
\email{shantanu.dave@univie.ac.at}
\address{\em{Department of Mathematics, University of Vienna, Austria}}
\thanks{}

\begin{abstract}
This article provides a survey of  abstract pseudo-differential operators based on \cite{Co-Mo}. In particular we  shall comment on the relationship between the algebraic and analytic concepts of order and under what conditions  they agree with each other. We shall also  note some regularity features that can be  provided in abstract generality and  correspond to geometrical  properties.
\end{abstract}
\maketitle
\section{Introduction}
The algebra of classical pseudo-differential operators on a closed manifold $M$  and its symbol calculus has been vastly useful in obtaining topological information about $M$ for example by means of index of some elliptic operator $D$ on it \cite{AS}.
  This inspired  some very abstract generalizations of both the concept of elliptic operators as well as the algebras of pseudo-differential operators which can be applied to the study of geometrical spaces other than smooth manifolds. In doing so  some deeper understanding of  concepts such as regularity has been obtained.  Most of these abstraction has been carried out in an operator-theoretic setup. Of particular interest, for our applications to certain singular spaces, are the abstract pseudo-differential operators of \cite{Co-Mo}.

Our aim here is to furnish an introduction to this abstract pseudo-differential calculus without the prerequisite of knowledge in noncommutative geometry.

\section{Notation}
In this article all algebras are over  the complex numbers $\CC$. A filtration on an algebra $\maD$ is an increasing sequence of subspaces, $\maD_0\subseteq \maD_1\subseteq \ldots$  such that
\[\maD=\bigcup_k \maD_k\quad  \textrm{and} \quad \maD_k\cdot \maD_l\subseteq \maD_{k+l}.\]
We then say that $\maD$ is a filtered algebra. Also if  an element $X$ is  a subspace $\maD_k$ we say that  $order(X)\leq k$.

For example the algebra of differential operators $\maD(M)$ on  any manifold $M$  can be given a filtration by  setting $\maD_k(M)$ to be the space of all operators  of order less than or equal to $k$.

We note that a symbol map on a filtered algebra $\maD$ is  a sequence of maps

$\sigma_k:\maD_k\rightarrow \maD_k/\maD_{k-1}.$

We say that an unbounded self-adjint operator $D$  on a Hilbert space  $\HH$, is positive if its spectrum $sp(D)\subseteq [0,\infty)$.
We shall use the convention  for  sign of the Laplace operator  on  a Riemannian manifold $M$ so that it  is a positive operator on $L^2(M)$. Thus for instace $\Delta_{\RR^n}=-\sum \frac{\partial^2}{\partial x_j^2}$.
\section{Abstract differential operators}
To motivate the notion of abstract differential operators we look at a familiar case.

First suppose that $M$ is a closed Riemannian  manifold and $\Delta$  the scalar Laplace operator. Then the algebra $\maD(M)$ of differential operator on $M$  and its filtration  with respect to the order of operators is  completely  determined by two pieces of data: a) the  algebra of smooth functions $\smooth(M)$  considered as operators by multiplication on $L^2(M)$ and, b) the  Laplace operator $\Delta$, as follows.
\begin{lem} The algebra $\maD(M)$ of differential operators on $M$ is the smallest  algebra of operators on $L^2(M)$  that contains $\smooth(M)$  and is closed under commutation with $\Delta$, that is
\[T\in \maD(M)\rightarrow [\Delta,T]\in\maD(M).\]
Furthermore the filtration on $\maD(M)$ is obtained by observing that the operators of order  $k$, the space $\maD_k(M)$, is generated by elements obtained by at most $k$  iterated commutators of a smooth function and $\Delta$.
\end{lem} 
 More generally, recall that  the differential operators $\maD(M:S)$ on sections of some Hermitian vector bundle $S$ on $M$,  can be generated by the space of bundle endomorphisms $End(S)$ and  all covariant derivations $\nabla_X$ with respect to some  connection $\nabla$ on $S$.  Here again if $\Delta$ is a Laplace type operator then the filtered algebra $\maD(M:S)$ is again completely determined by the operator $\Delta$ and the algebra of zero order operators  $\maD_0(M:S)=End(S)$ thanks to the Weitzenbock formula. 
 
 In \cite{Co-Mo} it is observed that  given such a description of the algebra of differential operators  using an operator $\Delta$  one can obtain an algebra of classical pseudodifferential operators. Our purpose is to explain this construction and some of its applications.
  
  We shall now abstract these ideas and define an abstract  notion of differential operators based on \cite{Co-Mo,Higson}. First let $\Delta $ be a positive self-adjoint operator on a Hilbert space $\HH$. Let $\HH^{\infty}$ be the intersection  $\cap_n  \operatorname{Domain} (\Delta^n)$. We shall assign $\Delta$ an order $r\geq0$ and for each real $m$ define the `Sobolev space' $H^m$ as the  completion of $Domain(\Delta^{\frac{m}{r}})$ with respect to the graph norm
  \[\|u\|_m^2:=\|u\|^2+\|\Delta^{\frac{m}{r}}u\|^2.\]
  Indeed in case of the concrete Laplace operator $\Delta$  this procedure defines the usual Sobolev norms  due to the  elliptic estimates, often called G{\aa}rding's inequalities.
  
 Back in our abstract setup one can define:
 \begin{definition}\label{DA}
A filtered algebra $\maD$  of operators on $\HH^{\infty}$ with filtration
\[\maD_0\subseteq \maD_1\subseteq \ldots \maD,\]
  is  an algebra of differential operators  for the operator $\Delta$  if the following hold.
  \begin{enumerate}
\item  If $T$ belongs to $\maD$ then so does $[\Delta,T]$ and
\[order([\Delta,T]\leq order(T)+r-1.\]
\item\label{elliptic} If $order(T)=k$ then $T$ extends to a bounded linear operator from $H^{s+k}$ to $ H^s$ for all $s$.
  \end{enumerate}
 \end{definition}

The condiition \ref{elliptic} above  is a reformulation of the  elliptic estimate of Laplace operator in the abstract settings. In general the algebra $\maD$  is not completely determined by commutators with the subalgebra of order zero operators $\maD_0$, but this would be the case of interest.

As is obvious from the definition  one can always set $\Delta$ to be a positive elliptic differential operator on some  manifold $M$. If $M$ is closed then, due to the elliptic estimates, the notion of Sobolev space remains unchanged.  In particular we recover the usual filtration iwth some possible rescaling. But it is a different matter  when the operator $\Delta$ is only hypoelliptic differential operator as seen  in the following simple geometrical example.
\begin{exm}\label{hopf} 
Consider the three sphere $S^3=\{z=(z_1,z_2)\in\CC^2|\,|z_1|^2+|z_2|^2=1\}$. Define a vector field 
\[V:=\left.\frac{d}{dt}(e^{it}z_1,e^{it}z_2)\right|_{t=0}.\]
Let $\Delta_{S^3}$ be the standard Laplacian on the  sphere $S^3$. Let us set $\Delta_{Hopf}:=\Delta_{S^3}+V^2-V^4$.
This new operator is clearly not elliptic. To see another form of it let us choose an orthonormal basis to the space $\{span(V(z)\}^{\perp}$ defining vector fields $X$ and $Y$  of unit norm, then
\[\Delta_{Hopf}=-X^2-Y^2-V^4.\]
Thus our new operator defines some weighted Sobolev spaces. Given  $u$ supported in a  small coordinate charts  along $X,Y,V$ can be  expressed by the norm:
\[\|u\|^2_m=\int |\hat{u}(\xi)|^2(1+\xi_1^2+\xi_2^2+\xi_3^4)^{\frac{m}{2}}d\xi,\]
where  $\xi_1,\xi_2,\xi_3$  are duel variables to $X,Y,V$ is the small neighbourhood.

The algebra of differential operators is the same as  the usual $\maD(S^3)$. Since the vector fields  $X,Y,V$ along with $\smooth(S^3)$ generate all differential operators, we  can get a new filtration  on   $\maD(S^3)$ as follows. We can assign the vector fields $X$  and $Y$  an order  $2$ and the vector field $V$  its usual order $1$.  Also assign the  multiplication by smooth functions an order of $0$.   The order of the operator $\Delta_{Hopf}$ in this new  world is still $4$    but it is elliptic with respect to the symbol calculus  obtained from this filtration.

This construction works in general for any foliated manifold which were source of inspiration in \cite{Co-Mo}.
\end{exm}

 Thus we see that the notions of filtration on the algebra of differential operators  naturally correspond to  notions of regularity in the Sobolev sense.
 
 In addition  the idea of  a symbol map is akin to the filtration used.  A (universal) symbol map on a filtered  algebra $\maA=\bigcup \maA_i~~\maA_0\subseteq \maA_1\subseteq \ldots $,  is just a  sequence of maps $\sigma+i: \maA_i\rightarrow \maA_i/\maA_{i-1}$. Thus by choosing a filtration we prescribe a choice of symbol map and even the meaning of ellipticity -which indicates invertibility modulo lower order operators.

 A different kind of symbol calculus was constructed by \cite{Getzler}. Consider $S$  the  spinor bundle over a spin manifold $M$.   Let us set $c:TM\rightarrow End(S)$ to be the Clliford action, $\nabla$ to be the spin connection and $D$ the Dirac operator.
The differential operators on $S$ can again be generated by the algebra $End(S)$ and the operator $\Delta:=D^2$ by iterated  commutators as in the above Lemma, but now we can choose from two different natural  filtrations on $\maD(M:S)$. The first is the filtration by the usual order of the differential operator. A more interesting filtration is described by Getzler \cite{Getzler} by exploiting the (representation theoretic) fact that 
\[End(S)=Clif(TM)\otimes End_{Cl}(S).\]
Here $Clif(TM)$ is the Cilfford algebra of $TM$ and its elements are prescribed an order $0$,  whereas  $End_{Cl}(S)$ are the Clifford module endomorphisms of $S$ and these elements are assigned an order $1$.   We can assign a filtration to $\maD(M:S)$ by  assigning an order $1$ to all covariant derivatives $\nabla_X$. We check that,
\[order([\Delta,T]\leq order(T)+1.\]
 This way of filtration give  both $\Delta=D^2$ and $D$  an order of ( less or equal to) $2$.  The resulting symbol calculus  provides a clearer  understanding of the topological index of $D$.
 
  We emphasize  again  that a different choice of filtration leads to  a different associated graded algebra $\mathfrak{G}r(\maD(M,S)$ and a different  symbol calculus.
  
  In the following section we shall start with an abstract filtered algebra of differential operators  associated to a positive operator $\Delta$ and prescribe it a calculus of pseudo-differential operators.

\section{Pseudo-differential operators}
As mentioned already we  shall  follow the Connes and Moscovici construction of pseudodifferential operators.

So given  an  algebra of  differential operators $\maD$ for a positive self-adjoint operator $\Delta$ in the sense of  Definition \ref{DA}  we shall first  extend the notion of order of an operator to a  bigger class of operators in  an obvious way. We say that an operator  $T$ on $\HH^{\infty}$ is of order $l$  if it extends to a bounded operator $H^{s+l}\rightarrow H^s$ for every $s\geq 0$.

\begin{definition}
We shall follow  \cite{Higson} in our formulation  and call an operator $T$ on $\HH^{\infty}$  a basic  pseudo-differential operator of order $k$ if for any $l\in\ZZ$  one can represent it in the form
\[T=X\Delta^{\frac{m}{r}}+R,\]
 where $X\in\maD$ and $order(X)\leq k-m$ and $order(R)\leq l$.
 
 More generally, a pseudo-differential operator is a finite linear combination of the basic pseudo-differential operators.
 \end{definition}
  At the first glance it is not at all clear that the space of pseudo-differential operators so defined is even an algebra. Hence we  introduce an asymptotic expansion for these abstract pseudo-differential algebras which is of  fundamental importance.
 
 \begin{definition}
 We say that an operator $T$  has an asymptotic expansion $T\approx \sum_j T_j$ if given any $l\in \ZZ$ there exists  an $N=N(l)$ such that for all $n\geq N$ the order relation
 \[order(T-\sum_{i\leq n}T_j)\leq l,\]
 is satisfied.
 \end{definition}

The basic non-commutative asymptotic expansion says that
\begin{prop}
Let $T$ be a pseudo-differential operator associated to $\maD$ and $\Delta$ as defined above. Let $\Delta_1=\Delta+K$ be an invertible operator such  that $K:\HH\rightarrow H^{\infty}$ is bounded\footnote{Such an opetator  $\Delta_1$ always exists by spectral thoery.}. Then
\[[\Delta_1^{-z},T]\approx \sum_j {-z \choose j} \nabla_1^j(T)\Delta_1^{-z-j}.\]

Here $\nabla_1(T):=[\Delta_1,T]$ and $\nabla_1^j$ is its $j-$th iterate.
\end{prop}
This proposition, proved  easily by means of holomorphic functional calculus, provides a powerful tool in understanding the abstract calculus of  pseudo-differential operators. This also shows in a standard way that the pseudo-differential operators form an algebra.

Let us go  back to our basic example when  $M$ is a  closed manifold  with $\Delta $ the Laplace operator and  $\maD(M)$  the honest to goodness  algebra of differential operators on $M$.  Recall now that  Kohn-Nirenberg \cite{KN} characterise a pseudo-differential operators of order $k$ as  an operator $T:\smooth(M)\rightarrow \smooth(M)$  such that for any  collection  of vector fields $X_1,X_2,\ldots X_k$,  the operator
\[[x_1,[X_2,[\ldots [X_k,T]\ldots ]]\]
extends to a bounded operator $H^{s+k}(M)\rightarrow H^s(M)$. This immediately tells us that our abstract construction does in fact give honest  pseudo-differential operators in the old framework  of closed manifold $M$, which is good news. The asymptotic expansion tells us that in fact we are talking about classical pseudo-differential operators.

\section{Generators for pseudo-differential operators}
So far we have  constructed our abstract algebra of pseudo-differential operators assuming that we are given two pieces of information, namely a filtered  algebra $\maD$ of abstract differential operators  associated with a positive self-adjoint operator  $\Delta$. In practise we  saw that  the subalgebra $\maD_0$ of order zero differential operators  is also a natural geometric entity like $ \smooth(M)$  or $End(S)$ which along side with $\Delta$  generates our differential operators $\maD(M)$. 

In general one can apply the following procedure. Fix a positive self-adjoint operator  $\Delta$  and the space of order zero  differential operators $\maD_0$. Then iteratively define the order $k$  differential operators by 
\begin{eqnarray}\label{gen_filtration}
\maD_k:=\sum \maD_i\cdot \maD_{k-i}+[\Delta,\maD_{k-1}]+\maD_0[\Delta,\maD_{k-1}]+\maD_{k-1}.
\end{eqnarray}
 In this situation we shall call the algebra $\maD=\bigcup_k\maD_k$ the filtered algebra generated by $\maD_0$ and $\Delta$.

This is all very good except that there is no reason to believe that such an algebraically generated space of differential operators would have the crucial  analytic property that each differential operator of order $k$ maps $H^{s+k}\rightarrow  H^s$ continuously.  To establish such a relationship between analysis  and algebraic-generation, we shall use the  language of spectral triples. 
 The setup of spectral triples envisioned by Connes provides some  deep insights.  A  principal example is provided by a Dirac operator on a spin manifold (see \cite{Connes} for more examples).

\begin{definition}
Suppose that we are given an algebra of bounded operators on $\HH$, and  an operator  $D$ on $\HH$ satisfying the following two conditions
\begin{enumerate}
\item  The  domain of $D$ is invariant under $\maA$, that is, $a\cdot \operatorname{Domain}(\Delta)\subseteq \operatorname{Domain}(\Delta)$ and the commutators $[a,D]$ extend to a  bounded operator  on $\HH$ for every $a\in\maA$.
\item  For every $a\in \maA$ the operator $a\cdot (1+D^2)^{-1}$ is a compact operator.
\end{enumerate}
Then we call the  structure $(\maA,\HH, D)$ a spectral triple.
\end{definition}

What interests us here  is that  a spectral triple provides  us with a candidate algebra of differential operators. 

 Now we can hopefully  choose the algebra $\maD_0(\maA):=\maA+[D,\maA]$ to be the algebra of order  $0$  differential operators  and hope that by picking $\Delta=D^2$  we can generate an appropriate  algebra of differential operators $\maD(\maA)$ using the procedure described by \eqref{gen_filtration}.
 
  So it turns out that the analytic condition  corresponding to the  elliptic estimates for $\Delta$ can be  reformulated as saying that all order zero differential operators  should also be order zero  pseudo-differential operators!

Remember that from  the Kohn-Nirenberg formulation the order zero pseudo-differential operators  on $M$ should satisfy the  condition that iterated commutators with vector fields should continue to be bounded operators on $L^2(M)$. 
Here such a condition is formulated as saying  that the space of order zero pseudo-differential operators  are all those bounded operators which extend to a bounded operator after finitely many  commutators with $|D|$. Thus we say that  $T\in B(\HH)$ is in $Op^0(\maA)$ if  any finite number of commutators,
\[[|D|,[|D|,[\ldots [|D|,T],\ldots  ]],\]
 extends to a  bounded operator on $\HH$. Then
\begin{definition}\label{regular_sp}
A  spectral triple is called regular if $\maD_0(\maA)\subset  Op^0(\maA)$.
\end{definition}

What this entails is the following:
\begin{prop}
A spectral triple defined by $\maA,\HH,D$ is regular if and only if the filtered algebra generated by $\maD_0=\maA+[D,\maA]$   with respect to $\Delta$ is an algebra of differential operators. 
\end{prop}
 This is essentially proved in Appendix B of \cite{Co-Mo} and is clearly mentioned in \cite{Higson}.
 
\section{Regularity}
So far we have seen that, given a positive self-adjoint operator $\Delta$ on a  Hilbert space $\HH$, we can construct an algebra of differential operators generated by a subalgebra $\maD_0$,  provided it corsponds to a regular spectral triple $(\maA,\HH, D)$.  We have also noted  that, with any  algebra of differential operator associated with $\Delta$ there is an algebra of   classical  pseudo-differential operators.   Lets us consider the  classical case  of the algebra $\maA$ equal $\smooth(M)$ where $M$ is a spin manifold. Here one can make another important observation namely that  the Frech\'et topology on $\smooth(M)$  can be obtained from the semi-norms $a\rightarrow \|\delta^n(a)\|$, where as before $\delta(a)=[|D|,a]$. This provides the differential structure of $M$ in the context of spectral triples. In fact Connes \cite{Connes} shows that  many other geometrical properties such as distance function on $M$ can easily be recovered  from the spectral triple data.

The question therefore is that, in cases where we do not have a manifold but something more singular albeit  with interesting geometrical properties like a length space (see for example \cite{Grant}), could we still apply a similar construction and obtain a calculus of pseudo-differential operator that encodes the regularity  in the structure concerned? This question is out of the scope of our present article and we shall focus on further understanding our abstract pseudo-differential operators for now.

Our first   desire is to check   if in this abstract settings the spectral zeta functions are able to produce interesting global invariants.  On a compact manifold, it is a non-trivial theorem due to Seeley \cite{Seeley} (and a different proof can be found in Guillemin \cite{Guillemin}) that the  complex powers $\Delta^{z}$ of an elliptic positive operator are in fact  pseudo-differential operators of order $\mathfrak{R}e(z)$ and hence when  $\mathfrak{R}e(z)$ is very large then
\[\zeta_A(z):=\operatorname{Tr}(A\Delta^{-\frac{z}{2}})\]
is an analytic function for any (classical) pseudo-differential operator $A$. Furthermore, it  has a meromorphic extension to whole of the complex plane with only simple poles. The assignment $\tau(A):= \underset{z=0}{\operatorname{res}}\zeta_A(z)$ provides a trace on the algebra of classical pseudo-differential operators and up to a constant multiple it is the unique trace on this algebra, referred to as  the noncommutative residue. The importance of noncommutative residue  lies in the fact that it is locally computable from the asymptotic  expansion of $A$ (even though the asymptotic expansion itself is not invariantly well-defined) and it provides global invariance as well.

The  main aim of \cite{Co-Mo} is to generalize this to a very abstract setting and provide under suitable conditions  a formula of a cyclic Chern-character of the class in  K-homology given by $\maA,\HH, D/|D|$  in terms of residues of certain zeta functions.

The fact that on an $n$ dimensional  manifold $M$   all pseudo-differential operators of order $<-n$ are of trace class,  whereas  the operators of order $\leq -n$ in general are in the  Dixmier ideal of $B(L^2(M))$, bears on the position of poles of the zetafunctions $\zeta_A(z)$.  Guillemin found a new proof of the Weyl's asymptotic formula, namely that the spectral counting function $N_D(\lambda)$, defined below, for  a positive elliptic differential operator $D$  is asymptotically given by:
\[N_D(\lambda):=\#\{\lambda_i\in sp(D)|\lambda_i\leq \lambda\}\simeq O(\lambda^{\frac{n}{order(D)}}),\]
by showing  that  the zeta function $\zeta_{Id}(z)=\operatorname{Tr}(D^{-z})$ has its highest pole at $z=\frac{n}{order(D)}$.

Of course, the Weyl asymtotics have practical applications including estimation of rates of convergences of numerical approximations.  It also provides a very clear understanding of regularity in the usual Sobolev sense.
It also provides a surprising  relationship between the poles of the analytically defined spectral zeta function and the geometric quantities such as dimension and volume. Such relationships  have been incarporated in the abstract setup of pseudo-differential operators. In practise, though it remains a difficult  proposition to prove a theorem like the theorem of Seeley \cite{Seeley, Guillemin} that is behind the analyticity of the zeta funciton on some half-plane of the complex domain, and even harder to establish a meromorphic continuation.


Thus we realize that  for a particular application of the abstract pseudo-differential calculus to a singular  space one would have to furnish statements that confirm that the axioms of regularity are satisfied. But in cases where the axioms can be ascertained, these abstract  pseudo-differential operators provide very useful geometrical information.

\section*{ Acknowledgements}
This work was supported by FWF grant Y237-N13 of the Austrian Science Fund. The author would like to thank Michael Kunzinger for his support and advice.

\bibliographystyle{gITR}


\end{document}